\documentclass[dvipsnames]{amsart}

\usepackage{tikz}
\usepackage[foot]{amsaddr}
\usetikzlibrary{plotmarks,
			  shapes.geometric,
			  positioning,
			  arrows,
			  decorations.pathreplacing,
			  decorations.markings,
			  automata}
\usepackage{placeins}
\usepackage{amsfonts}
\usepackage{mathtools}
\usepackage{hyperref}
\usepackage{stmaryrd}

\allowdisplaybreaks

\usepackage{amsmath}
\usepackage{amssymb}
\usepackage{algorithm}
\usepackage{algpseudocode}  
\usepackage{newtxtext,newtxmath}
\usepackage{textcomp}
\usepackage{graphicx}
\usepackage{bm}
\usepackage{multirow}
\usepackage{makecell}
\usepackage{booktabs}
\usepackage[caption=false]{subfig}
\usepackage[labelformat=parens,labelsep=quad,skip=3pt,justification=centering]{caption}
\usepackage{pgfplots} 
\usepackage{pgfplotstable} 
\pgfplotsset{compat=newest}
\pgfplotsset{plot coordinates/math parser=false}
\pgfplotsset{
tick label style={font=\footnotesize},
label style={font=\footnotesize},
legend style={font=\footnotesize}
}

\usepackage{latexsym}

\usetikzlibrary{plotmarks}
\usepackage{xcolor}
\definecolor{newcolor}{rgb}{.8,.349,.1}
\usepackage{tabularx}
\usepackage{tabu}
\usepackage{placeins}

\usetikzlibrary{
	plotmarks,
	decorations.pathreplacing,
	decorations.markings,
	automata,
	positioning,
	arrows,
	calc,
	spy,
	shapes,
	backgrounds,
	quotes,
	angles,
	arrows.meta
}
\usepgfplotslibrary{colormaps}
\pgfplotsset{
/pgfplots/colormap={jet}{rgb255(0cm)=(0,0,128) rgb255(1cm)=(0,0,255)
rgb255(3cm)=(0,255,255) rgb255(5cm)=(255,255,0) rgb255(7cm)=(255,0,0) rgb255(8cm)=(128,0,0)}
}
\pgfplotsset{%
  compat=newest,
  colormap={reversejet}{indices of colormap={
              \pgfplotscolormaplastindexof{jet},...,0 of jet}},
  colormap/orangeyellow/.style={
        colormap name=reversejet,
    },
}

\pgfplotsset{
/pgfplots/colormap={temp}{rgb255=(28.8,0,173.6) rgb255=(20,23.2,197.6) rgb255=(36.9,78.3,204)
rgb255=(48.8,108,204) rgb255=(69.6,140.8,204) rgb255=(93.6,168.8,204) rgb255=(122.4,188,204)
rgb255=(151.2,199.2,204) rgb255=(188,204,204) rgb255=(204,204,188) rgb255=(204,193.6,151.2)
rgb255=(204,171.2,122.4) rgb255=(204,137.6,93.6) rgb255=(204,96,52.2) rgb255=(204,48.4,48.8)
rgb255=(197.6,32,43.2) rgb255=(173.6,17.6,38.4) rgb255=(132.8,0,26.4)}
}

\usepackage{lipsum}
\usepackage{amsfonts}
\usepackage{graphicx}
\usepackage{epstopdf}

\ifpdf
  \DeclareGraphicsExtensions{.eps,.pdf,.png,.jpg}
\else
  \DeclareGraphicsExtensions{.eps}
\fi

\usepackage{amsopn}

\newtheorem{theorem}{Theorem}[section]

\theoremstyle{definition}

\theoremstyle{remark}
\newtheorem{remark}[theorem]{Remark}

\numberwithin{equation}{section}

\renewcommand{\bm}{\boldsymbol}

\begin{document}

\title[Binary mixture flow with free energy LBM]{Binary mixture flow with free energy lattice Boltzmann methods}

\author{Stephan Simonis$^{1,2,\ast}$}
\author{Johannes Nguyen$^{1,2}$}
\author{Samuel J. Avis$^3$}
\author{Willy Dörfler$^1$}
\author{Mathias J. Krause$^{1,2,4}$}

\email{stephan.simonis@kit.edu}
\email{j\_vietan@web.de}
\email{sam.avis@durham.ac.uk}
\email{willy.doerfler@kit.edu}
\email{mathias.krause@kit.edu}

\address{$^\ast$Corresponding author}
\address{$^1$Institute for Applied and Numerical Mathematics, Karlsruhe Institute of Technology, 76131 Karlsruhe, Germany}
\address{$^2$Lattice Boltzmann Research Group, Karlsruhe Institute of Technology, 76131 Karlsruhe, Germany}
\address{$^3$Department of Physics, Durham University, South Road, Durham DH1 3LE, United Kingdom}
\address{$^4$Institute of Mechanical Process Engineering and Mechanics, Karlsruhe Institute of Technology, 76131 Karlsruhe, Germany}


\date{{November 30, 2022}}


\begin{abstract}
We use free energy lattice Boltzmann methods (FRE LBM) to simulate shear and extensional flow of a binary mixture in two and three dimensions. 
To this end, two classical configurations are digitally twinned, namely a parallel-band device for binary shear flow and a four-roller apparatus for binary extensional flow. 
The FRE LBM and the test cases are implemented in the open-source C++ framework OpenLB and evaluated for several non-dimensional numbers. 
Characteristic deformations are captured, where breakup mechanisms occur for critical capillary regimes. 
Though the known mass leakage for small droplet-domain ratios is observed, suitable mesh sizes show good agreement to analytical predictions and reference results. 
\end{abstract}

\maketitle

\tableofcontents

\section{Introduction}
Mixture flows are omnipresent in nature and essential to many industrial processes. 
Taylor \cite{taylor1934formation} proposed machinery to examine the deformations of droplets induced by shear and extensional flow of binary multicomponent mixtures. 
The deformation is governed by the balance of outer forces and surface tension. 
Once this force balance is in favor of deformation, the droplet will break. 
Modifying the properties of the system, the breakup process can be adjusted in terms of number and size of resulting droplets. 
These phenomena are essential in manufacturing processes, for example in order to maximize the efficiency of creating emulsions \cite{calhoun2022systematic,wang2020modelling}.  
For the computer simulation of binary mixture flow several methods have been employed in the past. 
Due to the intrinsic parallelizability which enables the outsourcing of high performance computing (HPC) machinery, the lattice Boltzmann method (LBM) emerged as an unconventional alternative for multicomponent computational fluid dynamics (CFD). 
The popularity of LBM for CFD and beyond has increased significantly \cite{lallemand2021lattice}. 
Several data structures are available commercially and open-source. 
Exemplarily for the latter, the highly parallel C++ framework OpenLB \cite{krause2021openlb} has been successfully used for simulations of various transport processes on Top500 HPC machines (e.g. \cite{krause2021openlb,
haussmann2019direct,
simonis2020relaxation,
mink2021comprehensive,
dapelo2021lattice-boltzmann,
siodlaczek2021numerical,
simonis2021linear, 
haussmann2021fluid-structure,
simonis2022temporal, 
simonis2022constructing,
simonis2022forschungsnahe,
bukreev2022consistent}). 
Simulating multiphase and multicomponent flows in LBM is mostly based on a phase field model with diffuse interfaces. 
The interfacial zone brings forth additional physics captured by the Cahn\textendash Hilliard equation (CHE), though in turn upholds the high parallelizability of LBM. 
Several approaches for the underlying mixture dynamics exist \cite{huang2015multiphase}, for example the free energy model (FRE) \cite{swift1996lattice,semprebon2016ternary}. 
Tunable physical effects and top-down configuration of thermodynamics, are advantages of models akin to FRE. 
Albeit a high potential for numerical simulations, applications with FRE LBM for flows relevant to industrial processes are still rare. 

The dynamic effects on an immiscible binary component mixture can be abstracted into shear- and extension-dominated flows. 
For these two types of dynamic mixture flows, the present work implements and tests the FRE LBM with a simplistic binary fluid composition (equal density and viscosity). 
In particular, deformation as well as breakup phenomena are distinctively assessed to determine the models applicability for more complex applications. 
As such, we approve the suitability of the presented FRE LBM for numerically simulating these types of binary mixture flows via digitally twinning classical devices and comparing the results to references. 

This manuscript is structured as follows. 
Section~\ref{sec:methods} summarizes the methods, the numerical results are described in Section~\ref{sec:numerics} and Section~\ref{sec:conlusion} draws conclusions and closes the paper.

\section{Methodology}\label{sec:methods}

\subsection{Target equations}\label{subsec:teq}
A weakly compressible, isothermal fluid flow is described via the Navier\textendash Stokes equations (NSE) 
\begin{align}
	\partial_t \rho + \partial_{\alpha}(\rho u_\alpha) &= 0,   & \quad \text{in } \Omega \times I, \label{eq:nseCon}\\
	\partial_t(\rho u_\alpha) + \partial_{\beta}(\rho u_{\alpha} u _{\beta}) &= -\partial_{\beta} P^{\mathrm{th}}_{\alpha \beta} + \partial_{\beta} \nu (\rho \partial_{\alpha} u_{\beta} + \rho \partial_{\beta} u_{\alpha}), \quad & \text{in } \Omega \times I, \label{eq:nseMom}
\end{align}
where $\rho \colon \Omega \times I \to \mathbb{R}$ denotes the density, $\bm{u}\colon \Omega \times I \to \mathbb{R}^{d}$ is the fluid velocity, and $\nu > 0 $ is a kinematic viscosity, and \(\Omega \subseteq \mathbb{R}^{d}\) and \(I\subseteq \mathbb{R}^{+}\) denote space and time, respectively.
Here, $\mathbf{P}^{\mathrm{t h}} = \mathbf{P}^{\mathrm{c h e m}} + P \mathbf{I}_{d} \colon \Omega \times I \to \mathbb{R}^{d\times d}$ is the thermodynamic pressure tensor. 
For single phase and single component flow, $\mathbf{P}^{\mathrm{t h}}$ reduces to the isotropic pressure $P \mathbf{I}_{d}$. 
In case of a multicomponent flow the corresponding thermodynamics are introduced by the anisotropic chemical pressure tensor $\mathbf{P}^{\mathrm{c h e m}}$ \cite{kendon2001inertial}. 

To model a multicomponent system, capturing additional physics of the diffuse interface, the Cahn\textendash Hilliard equation (CHE) is introduced 
\begin{equation}
    \partial_t \phi + \bm{\nabla} \cdot (\phi \bm{u}) = M_\phi \bm{\Delta} \mu_\phi, \quad \text{in } \Omega \times I,
    \label{eq:che}
\end{equation}
where $\phi$ is the order parameter, $\mu_{\phi}$ denotes the chemical potential and $M_{\phi}$ is related to the mobility of the interface. 
Complemented with respective initial and boundary conditions, equations \eqref{eq:nseCon}, \eqref{eq:nseMom}, and \eqref{eq:che} are to be approximated with FRE LBM.

\subsection{Free energy lattice Boltzmann method}\label{subsec:free}
\begin{figure}[ht!]
  \centerline{
	\subfloat[\(D2Q9\)]{
		\includegraphics[]{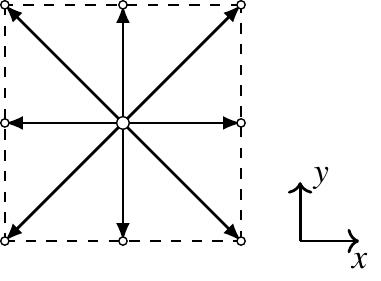}
    }
  }
  \centerline{
	\subfloat[\(D3Q19\)]{
 		\includegraphics[]{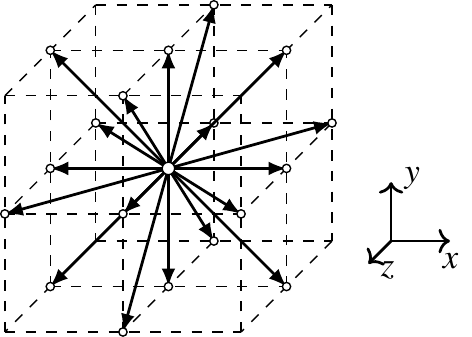}
	}
  }
  \caption{Discrete velocity sets.}
  \label{fig:discVeloSets}
\end{figure}
We assume a classical discretization of the phase space and the time domain \cite{kruger2017lattice} such that the following derivation is done completely in lattice units (\(\triangle t ^{\mathrm{L}} = 1 = \triangle x^{\mathrm{L}}\)). 

The present LBM is based on the lattice Boltzmann equation (LBE)
\begin{align}\label{eq:LBE}
f_{i} \left( \bm{x} + \bm{c}_{i} , t + 1 \right) = f_{i} \left( \bm{x}, t \right) + \Omega_{i}\left( \bm{x}, t\right) + S_{i} \left(\bm{x}, t\right), 
\end{align}
where \(i=0,1, \ldots, q-1\) counts the discrete velocities \(\bm{c}_{i}\) in \(DdQq\) and \(\bm{f}\left(\bm{x}, t\right) = ( f_{i}(\bm{x},t) )_{i}^{\mathrm{T}}\) denotes the populations in discrete space-time \((\bm{x}, t) \in (\Omega_{h}, I_{h}) \) and maps to \(\mathbb{R}^{q}\). 
The here used velocity stencils are given in Figure \ref{fig:discVeloSets}. 
The classical Bhatnagar\textendash Gross\textendash Krook (BGK) collision model \cite{bhatnagar1954model} reads 
\begin{align}
    \bm{\Omega}\left(\bm{x}, t\right) = -\frac{1}{\tau} \left[ \bm{f}\left(\bm{x}, t\right) - \bm{f}^{\mathrm{eq}}\left(\bm{x}, t\right)\right] , 
    \label{eq:bgk}
\end{align}
where $\tau > 0.5 $ denotes the relaxation time at which $\bm{f}$ relaxes towards the equilibrium  
\begin{align}
\label{eq:feq}
	f_i^{\mathrm{eq}} \left(\bm{x}, t\right)= w_i \rho \left(\bm{x}, t\right) \left[ 1 + \frac{c_{i \alpha} u_{\alpha}\left(\bm{x}, t\right) }{c_{\mathrm{s}}^2} + \frac{u_\alpha \left(\bm{x}, t\right) u_\beta \left(\bm{x}, t\right) \left( c_{i \alpha} c_{i \beta} - c_{\mathrm{s}}^2 \delta_{\alpha \beta} \right)  }{2c^4_s}\right]
\end{align}
with $w_i$ denoting the lattice weights and $c_{\mathrm{s}}$ is the lattice speed of sound. 
The term \(S_{i}\left(\bm{x}, t\right)\) obeys Guo's forcing scheme \cite{guo2002forcing} 
\begin{align}
    S_i \left( \bm{x}, t\right)= \left( 1 - \frac{1}{2\tau}\right) w_i \left[ \frac{\bm{c}_i - \bm{u} \left( \bm{x}, t\right)}{c_{\mathrm{s}}^2} + \frac{(\bm{c}_i \cdot \bm{u} \left( \bm{x}, t\right)) \bm{c}_i}{c_{\mathrm{s}}^4}\right] \cdot \bm{F} \left( \bm{x}, t\right), 
    \label{flowchart:Fi}
\end{align}
where \(\bm{F}\) is a force field. 
The macroscopic flow variables are recovered via the discrete moments of the populations 
\begin{align}
\rho \left(\bm{x}, t\right)& = \sum\limits_{i=0}^{q-1} f_{i}\left(\bm{x}, t\right), \label{eq:densMoment}\\
\rho \left(\bm{x}, t\right) \bm{u}\left(\bm{x}, t\right) & = \sum \limits_{i=0}^{q-1} \bm{c}_{i} f_{i}\left(\bm{x}, t\right) + \frac{1}{2}\bm{F} \left( \bm{x}, t\right), \label{eq:veloMoment}
\end{align}
respectively. 
As such, we use \eqref{eq:LBE} to approximate the conservative variables in the NSE \eqref{eq:nseCon} and \eqref{eq:nseMom}. 
Formal Chapman\textendash Enskog expansions are given in \cite{kruger2017lattice} and limit consistency is proven force-free in \cite{simonis2022limit}. 
Both establish the relation 
\begin{align}
\nu = c_{s}^{2}(\tau - 0.5) . 
\end{align}
The appropriate body force to recover the chemical pressure is defined as the residual 
\begin{align}
    F_{\alpha} &= -\partial_\beta (P^{\mathrm{th}}_{\alpha \beta} - c_{\mathrm{s}}^2 \rho \delta_{\alpha \beta}) \nonumber\\
    &= -\rho \partial_\alpha \mu_\rho - \phi \partial_\alpha \mu_\phi ,  
    \label{flowchart:F}
\end{align}
where \(\mu_{\rho}\) is specified below. 
For coupling the approximations of \eqref{eq:nseCon}, \eqref{eq:nseMom} and \eqref{eq:che} through the force \eqref{flowchart:F}, the thermodynamics can be consistently derived for a binary fluid mixture via the free energy functional \cite{kim2007phase,semprebon2016ternary}
\begin{align}\label{eq:functional}
        \Psi = \int\limits_{V}  \Bigl[ & \frac{\kappa_1}{32}(\rho + \phi)^2 (\rho + \phi -2)^2 + \frac{\alpha^2 \kappa_1}{8}(\bm{\nabla} \rho 
        + \bm{\nabla} \phi)^2  \nonumber \\ 
        & +  \frac{\kappa_2}{32}(\rho - \phi)^2 (\rho - \phi -2)^2 + \frac{\alpha^2 \kappa_2}{8}(\bm{\nabla} \rho - \bm{\nabla} \phi)^2 \Bigr] \mathrm{d}V,
\end{align}
where $\kappa_1$, $\kappa_2$, and \(\alpha\) are tunable parameters for the interface tension, and arguments are neglected where unambiguous. 
The auxiliary variables $\rho$ and $\phi$ are defined as
\begin{align}
\rho & = C_1 + C_2, \\
\phi & = C_1 - C_2, 
\end{align}
respectively, where $C_1$ and $C_2$ are the concentration fractions of the respective components.
The chemical potentials $\mu_\rho$ and $\mu_\phi$ are defined as functional derivatives of the free energy 
\begin{align}\label{flowchart:mu_rho}
        \mu_\rho = \frac{\delta \Psi}{\delta \rho} = & \frac{\kappa_1}{8}(\rho + \phi)(\rho + \phi - 2)(\rho + \phi - 1) \nonumber \\
        & + \frac{\kappa_2}{8}(\rho - \phi)(\rho - \phi - 2)(\rho - \phi - 1) 
               \nonumber \\  
               & +  \frac{\alpha^2}{4}\left[(\kappa_1 + \kappa_2)(-\bm{\nabla}^2 \rho) + (\kappa_2 - \kappa_1) \bm{\nabla}^2 \phi\right]
\end{align}
and
\begin{align} \label{flowchart:mu_phi}
        \mu_\phi = \frac{\delta \Psi}{\delta \phi} = & \frac{\kappa_1}{8}(\rho + \phi)(\rho + \phi - 2)(\rho + \phi - 1) \nonumber\\
        & - \frac{\kappa_2}{8}(\rho - \phi)(\rho - \phi - 2)(\rho - \phi - 1) \nonumber \\
                & + \frac{\alpha^2}{4}\left[(\kappa_1 + \kappa_2)(-\bm{\nabla}^2 \phi)  + (\kappa_2 - \kappa_1) \bm{\nabla}^2 \rho\right], 
\end{align}
respectively. 
The free energy of the system is minimized at equilibrium through the thermodynamic force induced by the chemical pressure. 
In case of a planar interface, the minimization of the free energy yields a simplified interface solution 
\begin{align}
\phi \left( x \right) = \mathrm{tanh} \left( \frac{x}{\sqrt{2}\xi} \right),
\end{align}
where $\xi$ is the interface width and the bulk components are identified by $\phi=\pm 1$ at $x=\pm \infty$ \cite{kusumaatmaja2010lattice}.
For approximating the CHE \eqref{eq:che}, a second population $g_i(\bm{x},t)$ is required such that its zeroth moment yields the order parameter
\begin{align} \label{eq:ordeMoment}
    \phi \left( \bm{x}, t\right)  = \sum\limits_i g_i \left( \bm{x}, t\right).
\end{align}
This population evolves according to a second LBE 
\begin{align}\label{eq:LBEg}
g_{i} \left( \bm{x} + \bm{c}_{i} , t + 1 \right) = g_{i} \left( \bm{x}, t \right) - \frac{1}{\tau_{g}} \left[ g_{i}\left(\bm{x}, t\right) - g_{i}^{\mathrm{eq}} \left(\bm{x}, t\right) \right], 
\end{align}
where \(\tau_g > 0.5\). 
To recover the coupled CHE \eqref{eq:che} in the continuum limit \cite{kruger2017lattice,semprebon2016ternary}, the corresponding equilibrium reads
\begin{align}\label{flowchart:geq}
	g_i^{\mathrm{eq}} \left( \bm{x}, t\right) = 
	w_i \Biggl[ \frac{\Gamma_{\phi} \mu_{\phi} \left( \bm{x}, t\right) }{c_{\mathrm{s}}^2} & + \frac{\phi\left( \bm{x}, t\right)  c_{i \alpha} u_\alpha\left( \bm{x}, t\right) }{c_{\mathrm{s}}^2} \nonumber \\
	& + \frac{\phi\left( \bm{x}, t\right)  u_\alpha\left( \bm{x}, t\right)  u_\beta \left( \bm{x}, t\right) \left (c_{i \alpha}c_{i \beta} - c^2_s \delta_{\alpha \beta})\right) }{2c^4_s} \Biggr],
\end{align}
if \( i\neq 0\), and 
\begin{align} 
	g_0^{\mathrm{eq}} \left( \bm{x}, t\right) = \phi\left( \bm{x}, t\right)  - \sum\limits_{i\neq 0}g_i^{\mathrm{eq}}\left( \bm{x}, t\right) , 	
\end{align} 
otherwise, where \(\Gamma_{\phi}\) relates to the mobility \cite{komrakova2014lattice} 
\begin{align}
M_{\phi} = \Gamma_{\phi} ( \tau_{g} - 0.5 ).
\end{align}
\begin{remark}
A flow with three or more components can be realized in a straightforward extension of the equation system \(\{\eqref{eq:nseCon}, \eqref{eq:nseMom}, \eqref{eq:che}\}\), by adding a similarly coupled CHE for each additional order parameter \cite{semprebon2016ternary} and thus one population for each additional component.  
\end{remark}

\subsection{Implementation}
The evolution equations for \(\bm{f} \) \eqref{eq:LBE} and \(\bm{g} = (g_{i})_{i}^{\mathrm{T}}\) \eqref{eq:LBEg} are split into local collision at \(\left(\bm{x}, t\right)\) which computes the post-collision populations \(\bm{f}^{\star}\) and \(\bm{g}^{\star}\), respectively, and streaming to evolve 
\begin{align}
\bm{f} ( \bm{x} + \bm{c}_{i}, t + 1 ) = \bm{f}^{\star}( \bm{x}, t ) 
\end{align}
and 
\begin{align}
\bm{g} ( \bm{x} + \bm{c}_{i}, t + 1 ) = \bm{g}^{\star} ( \bm{x}, t ) 
\end{align}
through the space-time cylinder.  
The computing steps implemented in one collision are summarized in Algorithm \ref{alg:freLBMkern} and the streaming is realized as a mere pointer shift \cite{kummerlander2022implicit}. 
The present work illustrates only the bulk solver. 

Standard LBM boundary methods are applied to impose velocity boundary conditions for the binary mixture. 
The macroscopic initial conditions are implemented via initialization of the populations to the corresponding equilibrium and alignment of kinetic moments through collisions \cite{mei2006consistent} preceding the actual simulation time horizon. 
\begin{algorithm}
\caption{FRE LBM bulk collision kernel}\label{alg:freLBMkern}
	\begin{algorithmic}[1]
	\Procedure{collide}{$\boldsymbol{f}$}  \Comment{Input: pre-collision $\boldsymbol{f}$ and \(\bm{g}\) at local node $(\boldsymbol{x},t)$} 
	\State compute zeroth moments: 
	\State \quad mixture density $ \rho \gets \sum_{i} f_{i} $ \Comment{\eqref{eq:densMoment}}
	\State \quad order parameter $ \phi \gets \sum_{i} g_{i} $\Comment{\eqref{eq:ordeMoment}}
	\State compute potentials: $ \mu_{\rho} \gets ( \rho, \phi ) $ and $ \mu_{\phi} \gets ( \rho, \phi ) $ \Comment{\eqref{flowchart:mu_rho}, \eqref{flowchart:mu_phi}}
	\State Guo forcing, via: 
	\State \quad force $ \bm{F} \gets ( \rho, \mu_{\rho}, \phi, \mu_{\phi} ) $ \Comment{\eqref{flowchart:F}}
	\State \quad velocity $ \bm{u} \gets ( \bm{f}, \rho, \bm{F} ) $ \Comment{\eqref{eq:veloMoment}}
	\State \quad force term $\bm{S} \gets  ( \bm{u}, \bm{F} ) $
	\Comment{\eqref{flowchart:Fi}}
	\State compute equilibria: $\bm{f}^{\mathrm{eq}} \gets  ( \rho, \bm{u}  ) $ and $\bm{g}^{\mathrm{eq}} \gets ( \phi, \mu_{\phi}, \bm{u} )$ \Comment{\eqref{eq:feq} , \eqref{flowchart:geq}}
	\State local collision of $\bm{f}^{\star} \gets ( \bm{f}, \bm{f}^{\mathrm{eq}}, \bm{S} )$ and $\bm{g}^{\star} \gets ( \bm{g}, \bm{g}^{\mathrm{eq}} )$ \Comment{\eqref{eq:LBE}, \eqref{eq:LBEg}}
	\State \textbf{return} $\bm{f}^{\star}$, $\bm{g}^{\star}$ 
	\EndProcedure \Comment{Output: post-collision $\bm{f}^{\star}$ and \(\bm{g}^{\star}\) at local node $( \boldsymbol{x}, t )$ }
	\end{algorithmic}
\end{algorithm}

\section{Numerical results}\label{sec:numerics}
\begin{figure}[ht!]
	\centering
	\subfloat[Shear flow]{
		\includegraphics[]{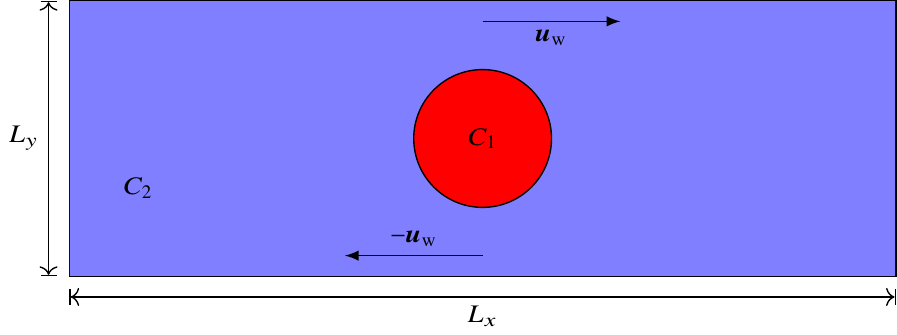}
	} \\
	\subfloat[Extensional flow]{
		\includegraphics[]{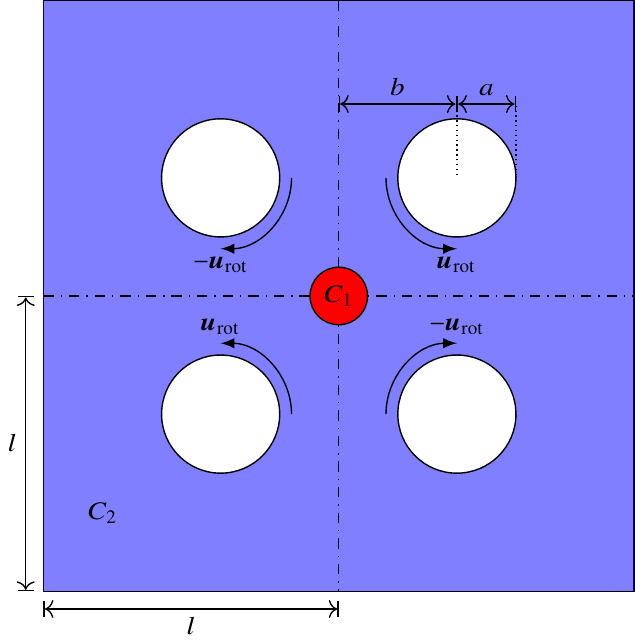}
	}
\caption{Geometric setup of numerical test cases for binary flow in two dimensions. 
Scales differ for the purpose of representation.}
\label{fig:apparatus}
\end{figure}

To assess the capability of the FRE LBM for recovering shear and extensional binary flows, we emulate Taylor's parallel-band and four-roller devices \cite{taylor1934formation} by means of numerical simulations in two dimensions (2D). 
The former is also tested for three dimensions (3D). 
The geometric setup of the 2D simulations is sketched in Figure~\ref{fig:apparatus}. 

All computations are done with OpenLB release 1.4 \cite{krause2020openlb14} on several HPC machines, either using up to 16 nodes with five quad-core Intel Xeon E5-2609v2 cores each, or a maximum of 75 nodes with respectively two deca-core Intel Xeon E5-2660v3. 

The deformation of the \(C_{1}\) droplet 
\begin{align}
D=\frac{L-B}{L+B}, 
\end{align}
where \(L\) is the longer axis halved and \(B\) the shorter one, is measured via intrinsic functors of OpenLB \cite{krause2021openlb}. 
In case of a horizontally measured inclination angle \(\theta = 0^{\circ}\), an interpolation along the space directions recovers the location of the interface. 
If a simultaneous deformation and inclination of the droplet occurs, \(L\) and \(B\) are approximated through concentric circles and \(\theta\) is computed at the intersection point. 

Though essential differences between 2D and 3D deformations are known \cite{soligo2020deformation}, certain non-dimensional regimes still allow a side-by-side comparison. 
Based on that we compare the FRE LBM solutions to 3D reference computations \cite{komrakova2014lattice,li2000numerical} and analytical predictions \cite{shapira1990low,taylor1934formation}. 

\subsection{Binary shear flow}\label{subsec:binaryShear}
We define the non-dimensional Reynolds number, capillary number, P{\'e}clet number, and Cahn number, respectively as 
\begin{align}
R\!e &= \frac{\gamma a^{2}}{\nu}, \\
C\!a &= \frac{a \gamma \mu_{\mathrm{c}} }{ \sigma}, \\
P\!e &= \frac{\gamma a \xi}{M_{\phi} A}, \\
C\!h &= \frac{\xi}{a},
\end{align}
where 
\(\gamma\), 
\(a\), 
\(\mu\), 
\(\sigma\), 
\(\xi\), and 
\(A\) 
are 
shear rate, 
droplet radius, 
viscosity, 
surface tension, 
interface thickness and a mobility parameter, respectively. 
The ratios of viscosity and density of the components are unity.  
\begin{figure}[ht!] 
\centerline{
	\includegraphics[]{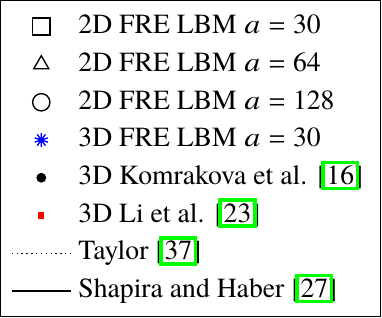}
}
\centerline{
	\subfloat[Deformation]{
		\includegraphics[]{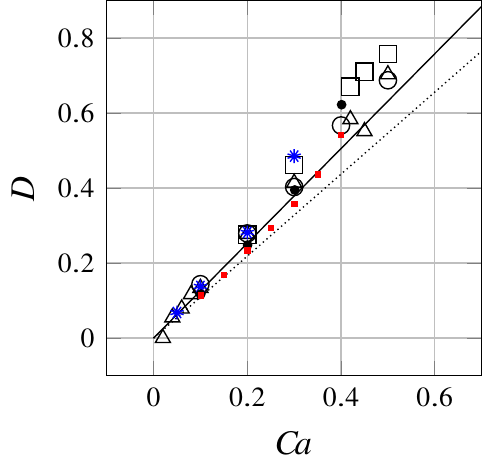}
		\label{fig:9Ca_a}	
	}  \hspace{2em}  
	\subfloat[Inclination]{
		\includegraphics[]{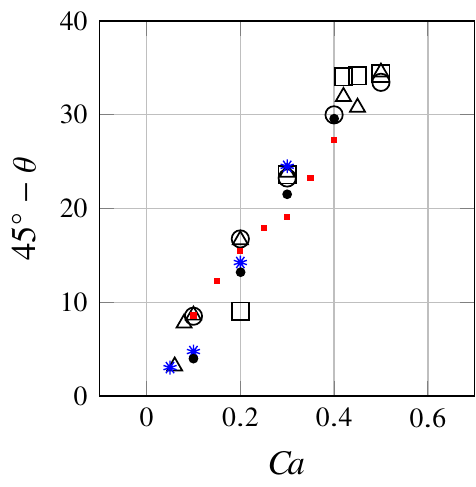}
		\label{fig:bsfSteadyState}
	} 
}
\caption{Deformation and inclination of a droplet in binary shear flow simulated with FRE LBM for varying capillary numbers.}
\end{figure}

\subsubsection{Steady state validation} 
In the case of \(C\!a < C\!a_{\mathrm{c}}\), the droplet deforms until it reaches a steady state. 
For \(R\!e = 0.1\), \(P\!e = 0.43\), and \(C\!h = 0.0379\), \(C\!a\) is varied over the interval \((0.02,0.6)\). 
The results are plotted in Figure~\ref{fig:bsfSteadyState} and agree well with the literature for small \(C\!a\). 
Refining the mesh over several droplet radii in lattice units \(a=30, 64, 128\), indicates convergence to the reference results and allows to simulate validly for higher \(C\!a\). 
\begin{figure}[ht!]
   \centerline{
	 \subfloat[\centering $\overline{t}=8$]{{\includegraphics[width=0.95\textwidth]{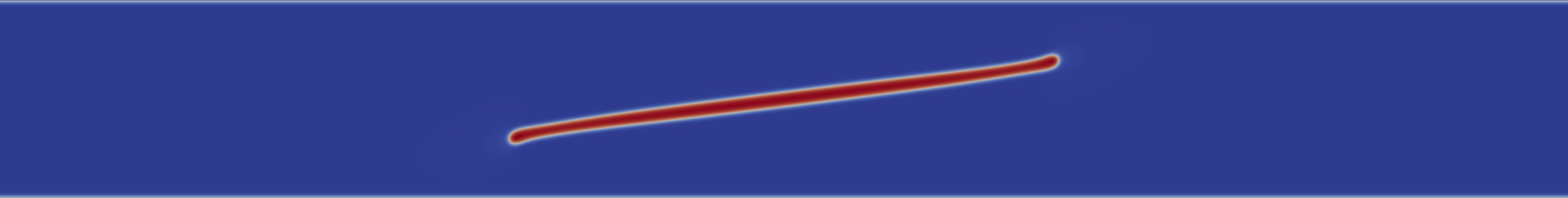} }}}
	\centerline{
   	 \subfloat[\centering $\overline{t}=14$]{{\includegraphics[width=0.95\textwidth]{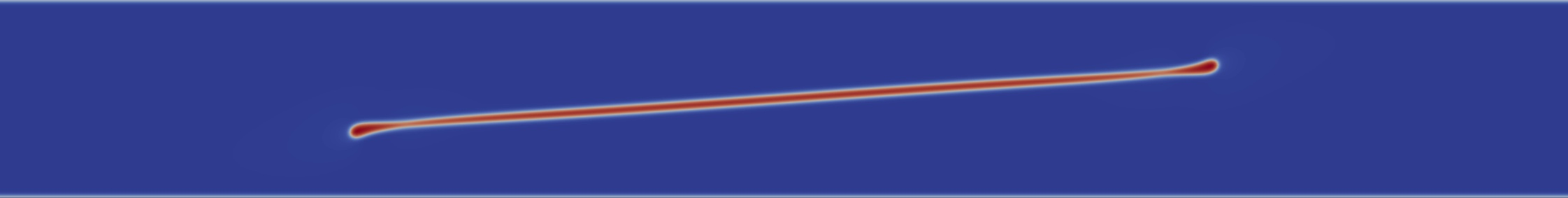} }}}
   \centerline{
	 \subfloat[\centering $\overline{t}=15$]{{\includegraphics[width=0.95\textwidth]{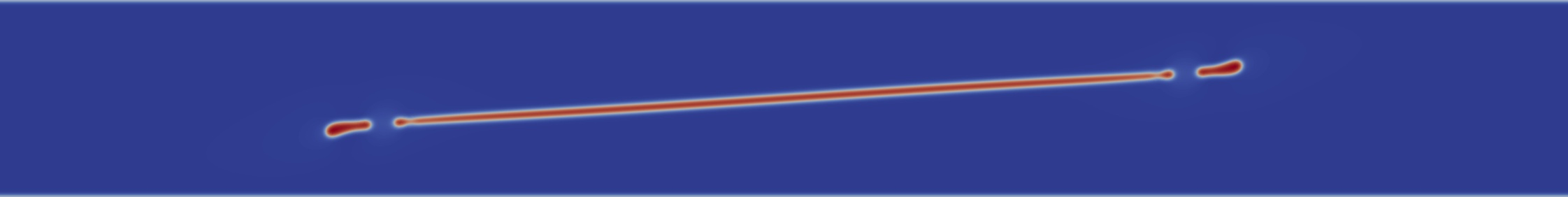} }}}
   \centerline{
	 \subfloat[\centering $\overline{t}=18.8$]{{\includegraphics[width=0.95\textwidth]{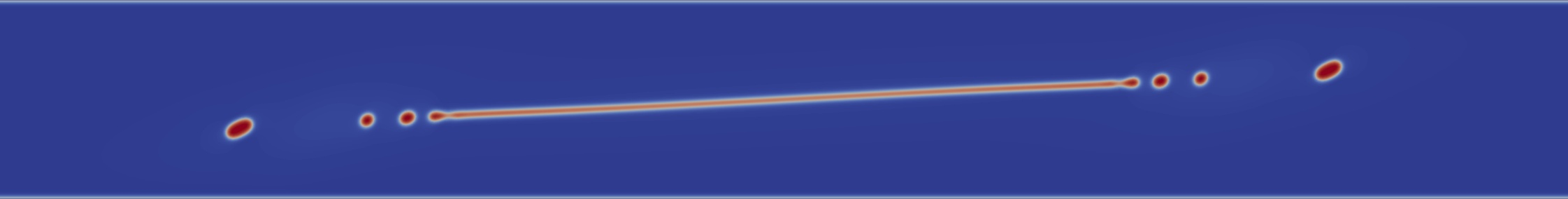} }}}
   \centerline{
	 \subfloat[\centering $\overline{t}=19.4$]{{\includegraphics[width=0.95\textwidth]{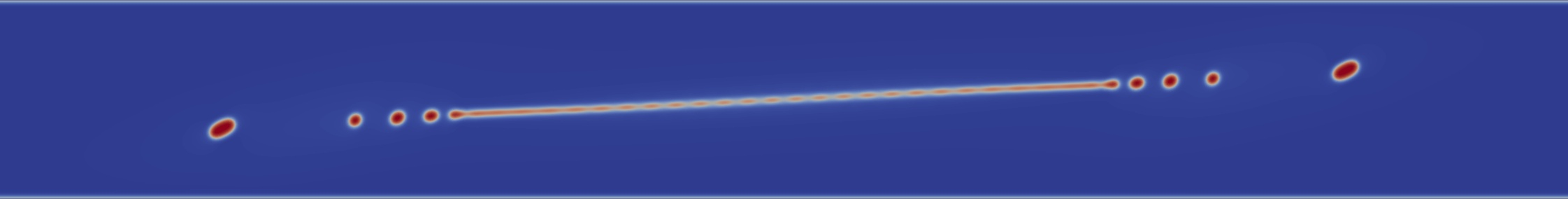} }}}
   \centerline{ 
   \subfloat[\centering $\overline{t}=20.4$]{{\includegraphics[width=0.95\textwidth]{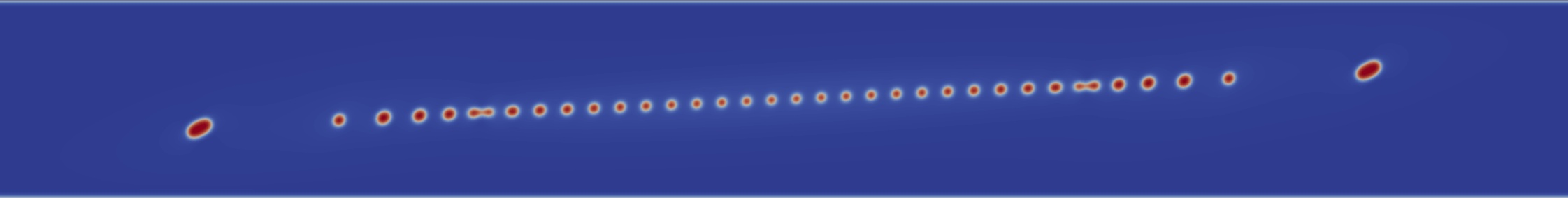} }}}
    \caption{Droplet breakup in 2D binary shear flow at $C\!a = 3.5$. 
    Components \(C_{1}\) (red) and \(C_{2}\) (blue) are plotted at normalized time steps.}
    \label{fig:capillaryWave}
\end{figure}

\begin{figure}[ht!] 
   \centerline{
	 \subfloat[\centering $\overline{t}=3$]{{\includegraphics[width=0.95\textwidth]{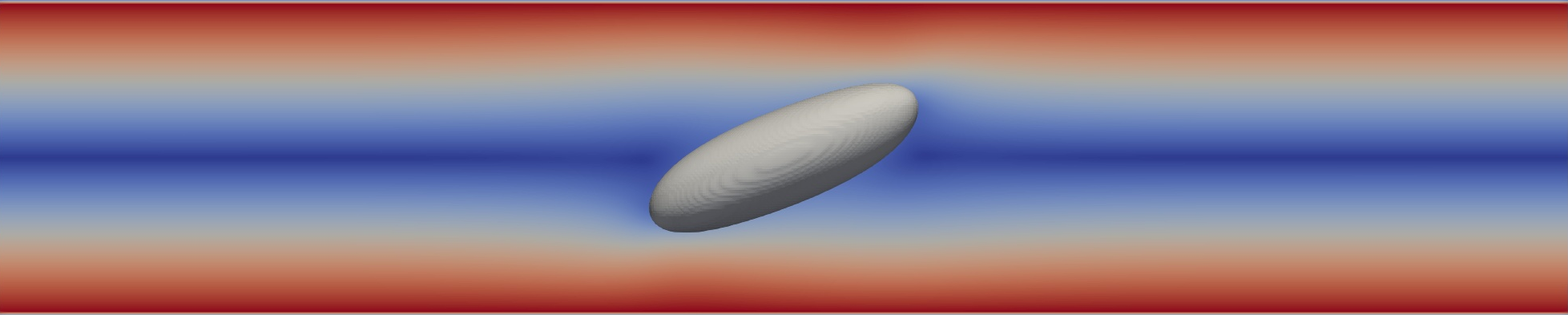} }}}
   \centerline{
	 \subfloat[\centering $\overline{t}=26.25$]{{\includegraphics[width=0.95\textwidth]{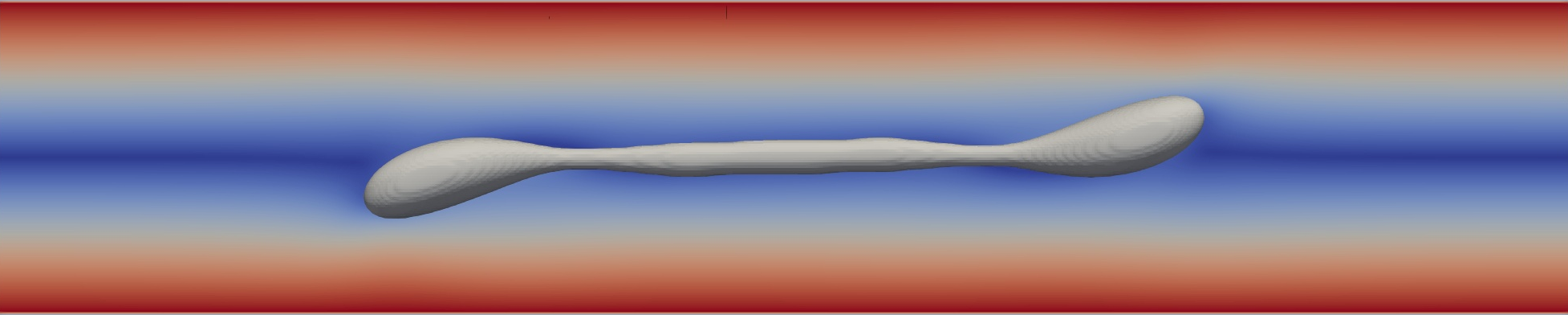} }}}
	    \centerline{
	 \subfloat[\centering $\overline{t}=27.5$]{{\includegraphics[width=0.95\textwidth]{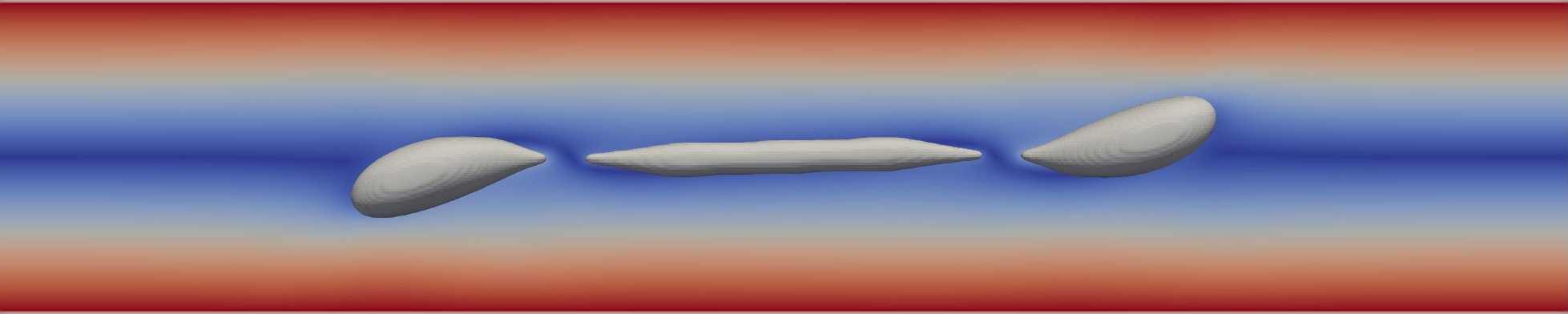} }}}
	\centerline{
	 \subfloat[\centering $\overline{t}=30$]{{\includegraphics[width=0.95\textwidth]{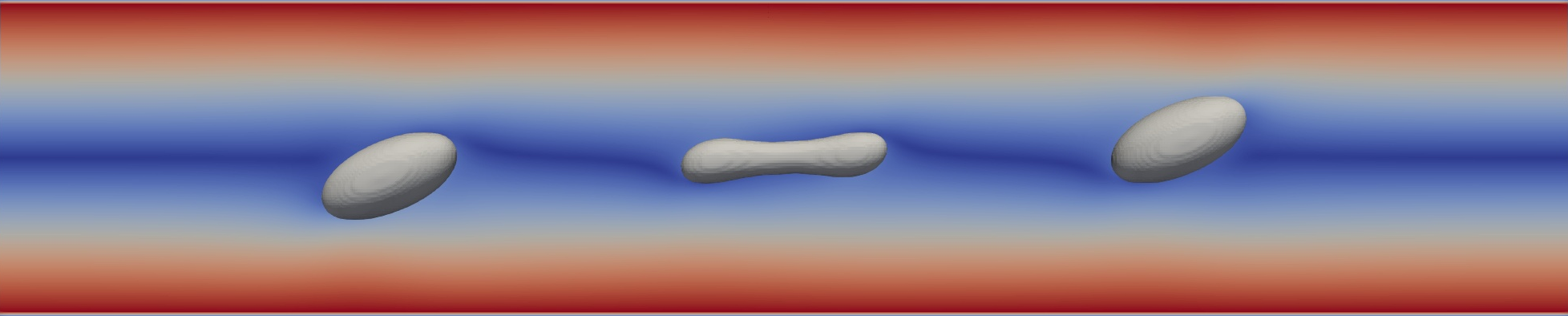} }}}
   \centerline{
	 \subfloat[\centering $\overline{t}=32.25$]{{\includegraphics[width=0.95\textwidth]{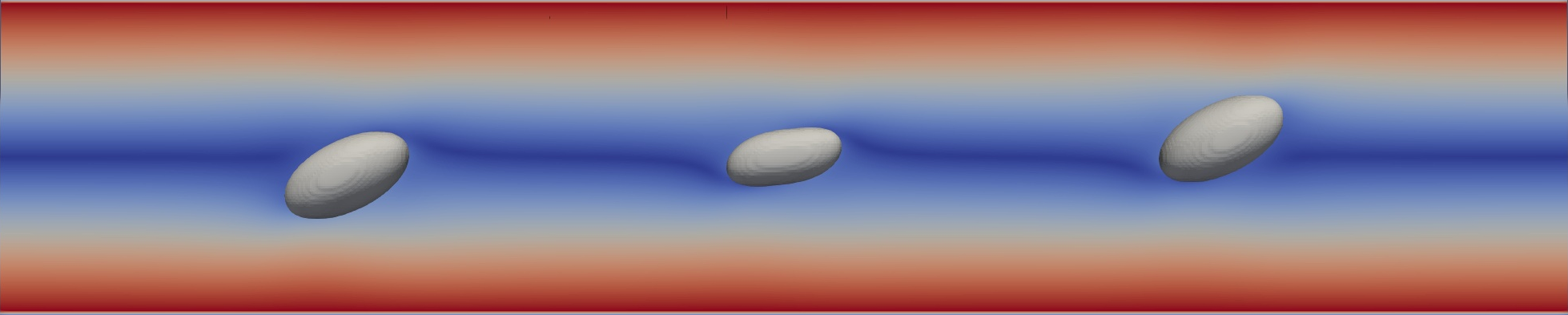} }}}
  	\centerline{
		\includegraphics[]{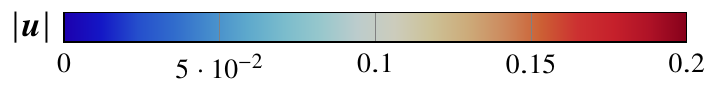} 
	} 
    \caption{Droplet breakup in 3D binary shear flow at normalized time steps for \(R\!e=0.0625\), \(C\!a = 0.42\), \(C\!h = 0.0379\), \(P\!e = 0.43\).}%
	\label{fig:breakup}
\end{figure}

\subsubsection{Breakup}
The breakup occurs for \(a=30\) approximately at \(C\!a_{\mathrm{c}}\sim 0.7\) (2D) and \(C\!a_{\mathrm{c}} \sim 0.42\) (3D).  
Three different categories exist \cite{bigio1998predicting,zhao2007drop}, namely
\begin{enumerate}
\item[(i)] pseudo steady-state $C\!a \sim C\!a_c$,
\item[(ii)] end-pinching $ C\!a_c < C\!a < 2 C\!a_c$, and
\item[(iii)] capillary wave breakup $C\!a > 2 C\!a_c$. 
\end{enumerate}
The bounds between these regimes however are not sharp, such that the droplet may pass through multiple types during the breakup process. 
Due to differences between 2D and 3D droplet deformations, \(C\!a_{\mathrm{c}}\) in 2D is significantly larger than in 3D. 
The latter agrees well with the literature \cite{komrakova2014lattice}, and so does the breakup scenario (see Figure~\ref{fig:breakup}). 
In 2D for \(C\!a=5C\!a_{\mathrm{c}}\) at $R\!e = 1$, $P\!e=0.2$ and $a=40$ we observe end-pinching as well as a capillary wave breakup (see Figure \ref{fig:capillaryWave}).

\subsection{Binary extensional flow}\label{subsec:fourRoller}
The sizing of the four-roller device ensures a uniform extension rate \(\epsilon\) \cite{higdon1993kinematics} which now replaces \(\gamma\) \cite{park2019taylor,
tretheway2001deformation,bentley1986experimental} in the \(\pi\)-group. 
For fixed $R\!e=0.0625$, $C\!h=0.57$, $P\!e=0.1$ and $C\!a \in [0.01, 0.3]$ the droplet is observed to break for $C\!a > 0.25$.
\begin{figure}
\centerline{
		\includegraphics[]{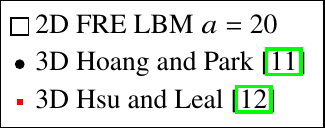} 
}\vspace{1em}  
\centerline{ 
		\includegraphics[]{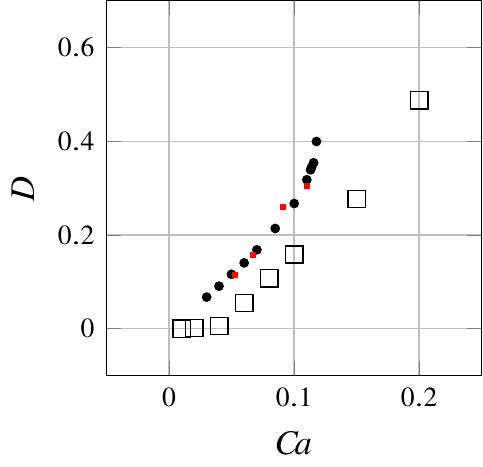} 	
} 
\caption{Deformation of a droplet in binary extensional flow simulated with FRE LBM for varying capillary numbers.}
\label{fig:frm_def}
\end{figure}

\subsubsection{Steady state validation}
The droplet radius is set to $a=20$, which corresponds to a ratio of \(40\) between domain length and radius.
Figure \ref{fig:frm_def} summarizes the deformation in the subcritical capillary regime.
For $C\!a=0.01, 0.02, 0.04$ the droplet shows little to no deformation. 
Beginning at $C\!a=0.05$ the deformation becomes significant and increases rapidly with increasing \(C\!a\) and with a considerably faster rate than in the shear flow. 
Our simulation results and the reference data from \cite{park2019taylor,hsu2009deformation} agree from the perspective of an overall trend but differ at individual values. 
Based on the same reasoning as above, we conclude that a 3D extensional flow produces a higher deformation at lower \(C\!a\) than in 2D. 
\begin{figure}[ht!]
   \centerline{
	 \subfloat[\centering $\overline{t}=0.25$]{{\includegraphics[width=.4\textwidth]{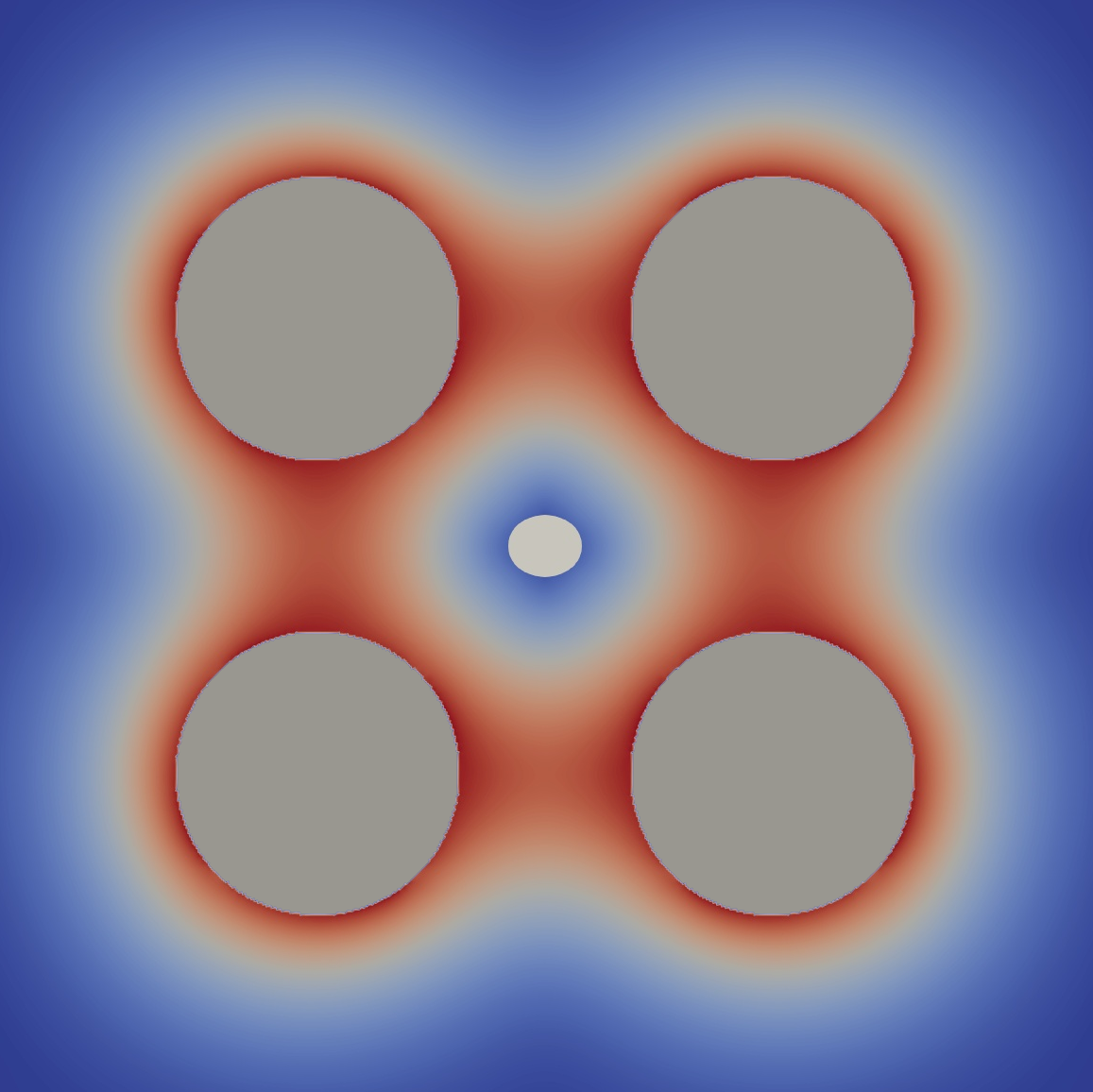} }}
	 \subfloat[\centering $\overline{t}=1.25$]{{\includegraphics[width=.4\textwidth]{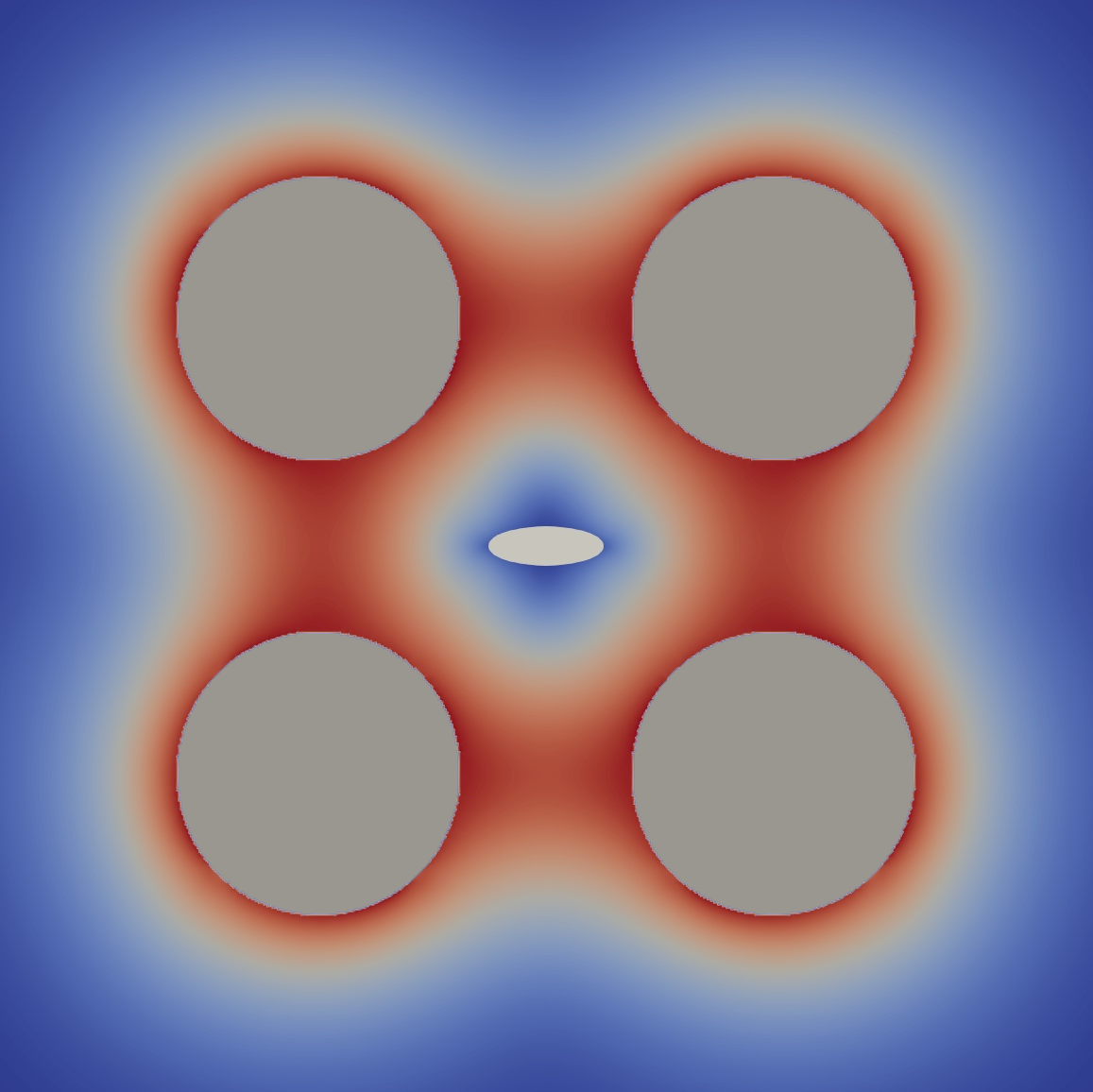} }}}
   \centerline{
	 \subfloat[\centering $\overline{t}=4$]{{\includegraphics[width=.4\textwidth]{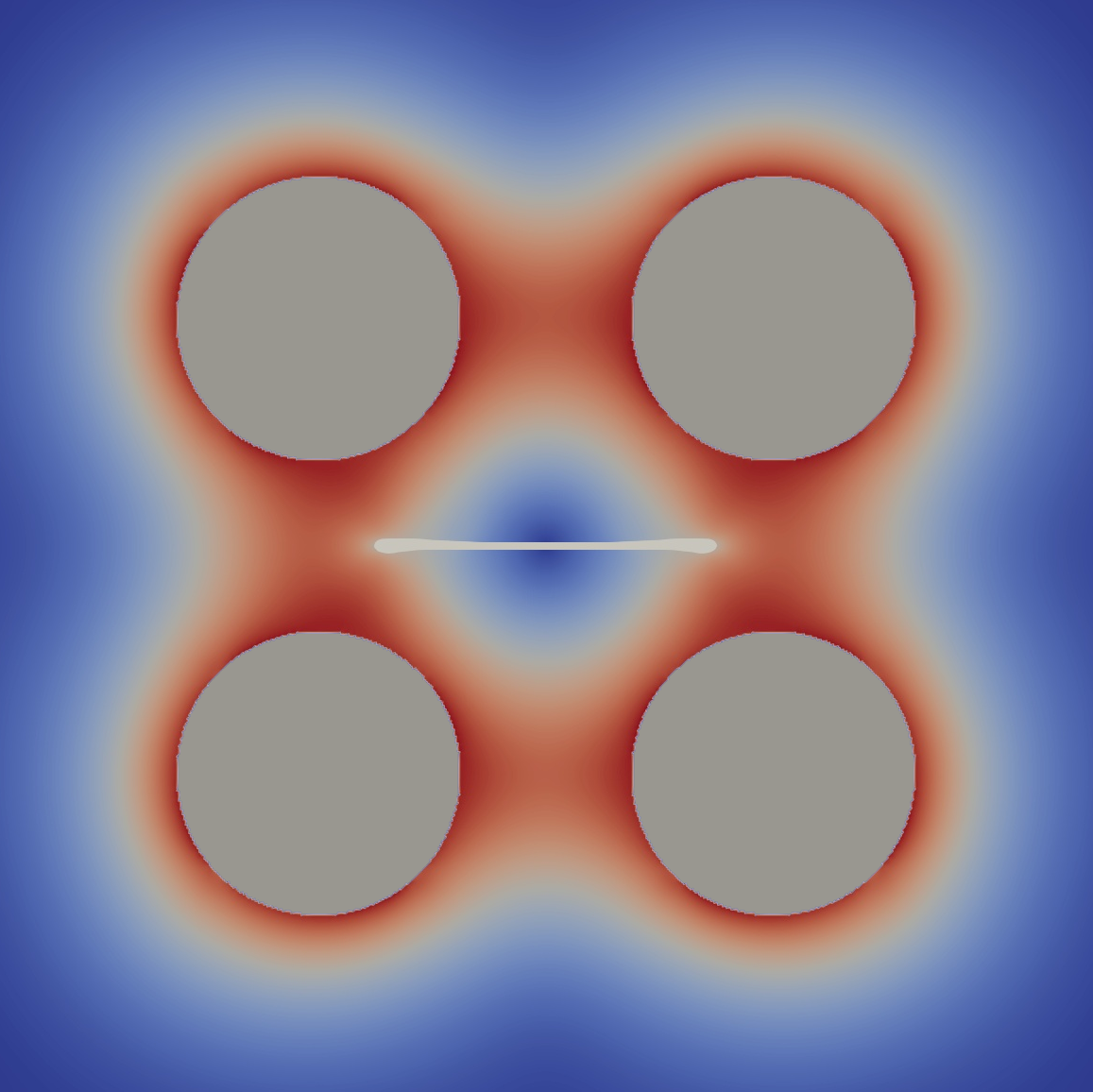} }}
	 \subfloat[\centering $\overline{t}=4.75$]{{\includegraphics[width=.4\textwidth]{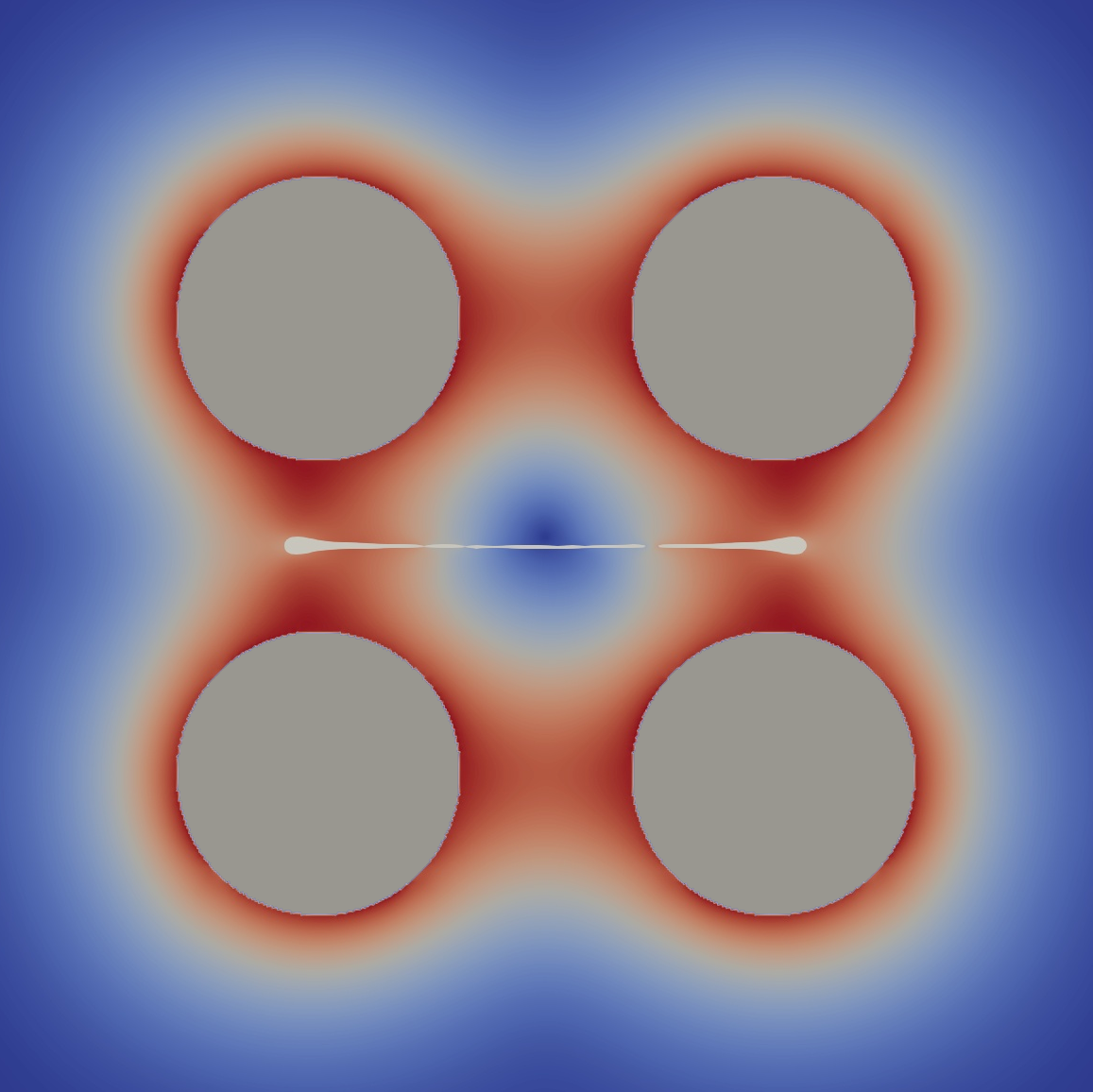} }}  }
   \centerline{
	 \subfloat[\centering $\overline{t}=5$]{{\includegraphics[width=.4\textwidth]{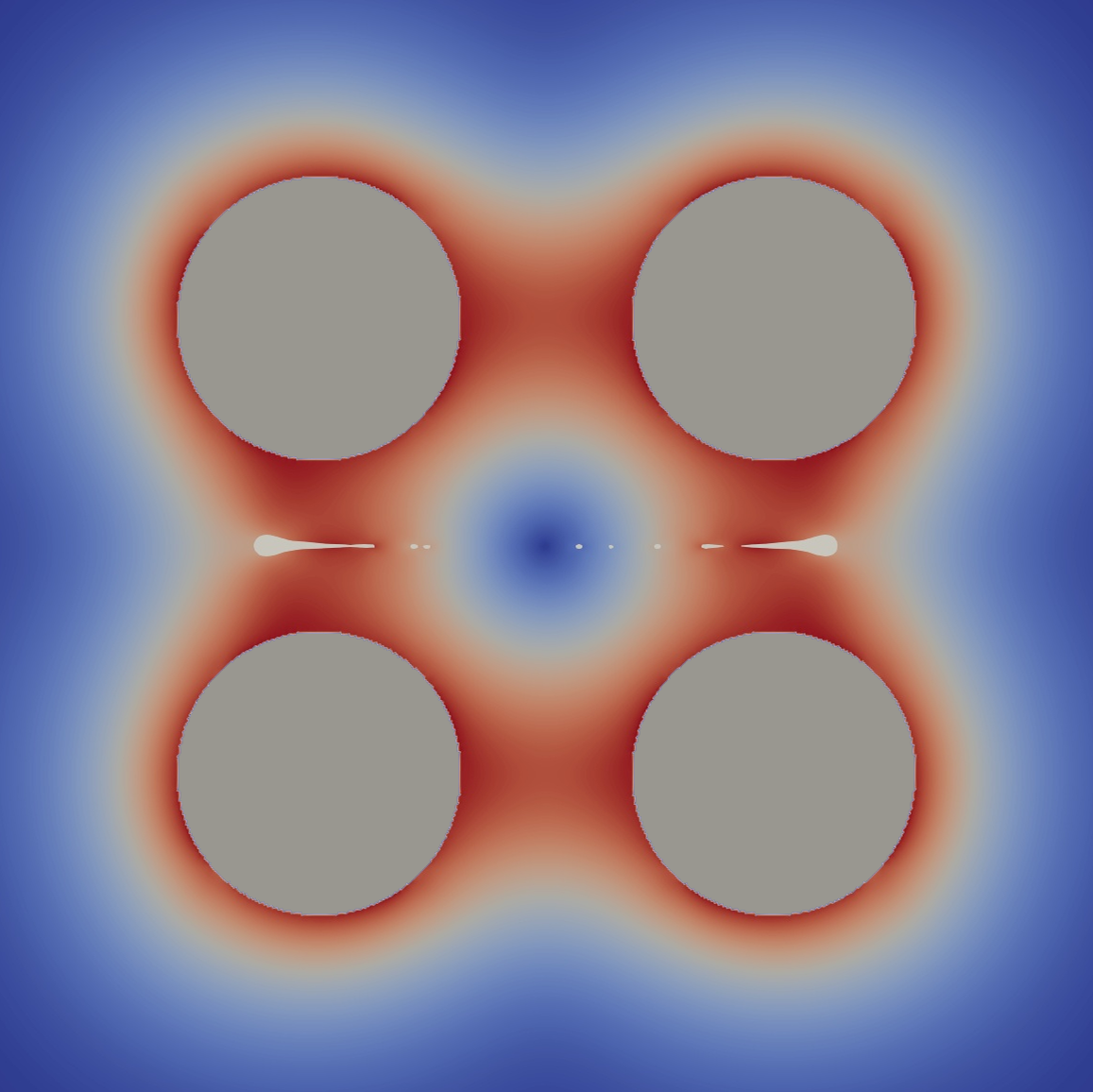} }}     
	 \subfloat[\centering $\overline{t}=5.5$]{{\includegraphics[width=.4\textwidth]{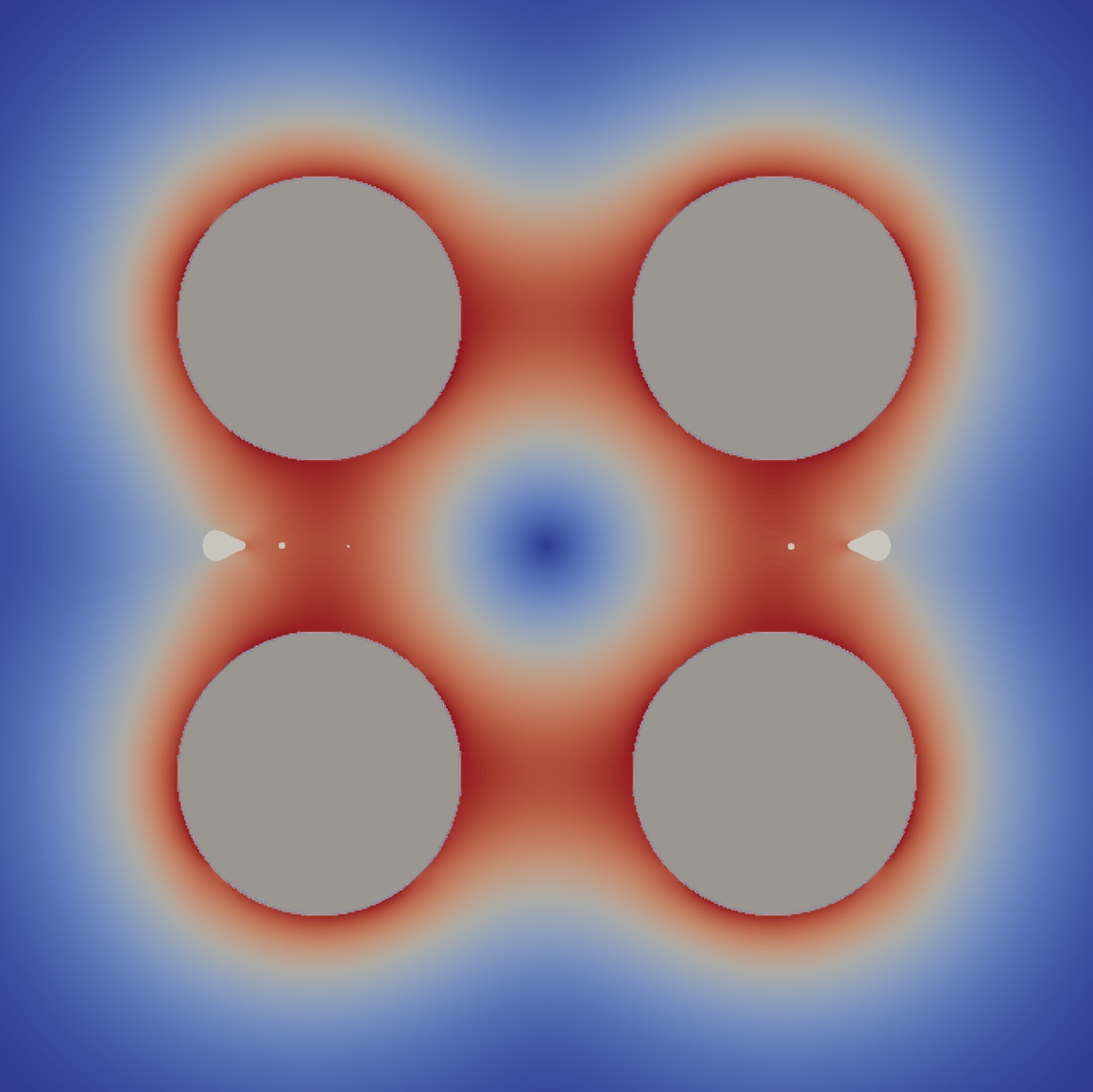} }} 
	 }
	  \centerline{
		\vspace{0.4em}
    	\includegraphics[]{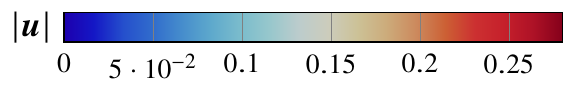} 
	}
    \caption{Droplet breakup in 2D binary extensional flow at normalized time steps for \(R\!e=0.0625\), \(C\!a = 0.42\), \(C\!h = 0.057\), \(P\!e = 0.43\).}%
    \label{fig:frmBreakup}
\end{figure}

\subsubsection{Breakup} 
Figure \ref{fig:frmBreakup} visualizes a breakup for $C\!a = 0.42$. 
At first, the droplet stretches into a long thread of equal width with rounded ends. 
Instead of end-pinching, the droplet breaks by overstretching, forming several sub-satellite droplets. 
After the droplets exceed the rollers' gap, their velocity declines and the tails retract quickly. 
Despite the successful reproduction of these phenomena, the results show heavy satellite shrinkage due to the high domain-droplet ratio \cite{zheng2014shrinkage}.

\section{Conclusion}\label{sec:conlusion}
We set up a FRE LBM algorithm implementation in OpenLB and test it for shear and extensional binary mixture flow in two and three dimensions. 
Taylor's parallel-band and four-roller devices are digitally twinned in a simplified manner and used for validation of the numerical scheme. 
To the authors' knowledge, the present work is the first application of LBM for simulating a four-roller apparatus. 
Characteristic deformations for steady states and breakup scenarios in critical capillary regimes are captured. 
Though the known satellite loss for very small droplet-domain ratios is observed, with suitably fine meshes we find good agreement to references. 

To push forward the FRE LBM application in industrial processes, future studies should include model extensions to non-uniform viscosity and density ratios as well as conservative discretizations.

\section*{Acknowledgements}
The authors thank Halim Kusumaatmaja for valuable discussions. 
S.S. thanks Taehun Lee for pointing out valuable references. 
The authors acknowledge support by the state of Baden-Württemberg through bwHPC. 
Parts of this work were performed on the supercomputer ForHLR II funded by the Ministry of Science, Research and the Arts Baden-Württemberg and by the Federal Ministry of Education and Research.

\section*{Author contributions}
Conceptualization: \textbf{S.S.}; 
Methodology: \textbf{S.S.}, \textbf{J.N.}, \textbf{S.J.A.}; 
Software: \textbf{S.S.}, \textbf{J.N.}, \textbf{S.J.A.}, \textbf{M.J.K.};  
Validation: \textbf{S.S.}, \textbf{J.N.};  
Formal analysis: \textbf{S.S.}, \textbf{J.N.};  
Investigation: \textbf{S.S.}, \textbf{J.N.};  
Data Curation: \textbf{S.S.}, \textbf{J.N.}; 
Writing - Original draft: \textbf{S.S.}; 
Visualization: \textbf{S.S.}, \textbf{J.N.};  
Supervision: \textbf{S.S.};  
Project administration: \textbf{S.S.};  
Writing - Review {\&} Editing: \textbf{S.S.}, \textbf{J.N.}, \textbf{S.J.A.}, \textbf{W.D.}, \textbf{M.J.K.}; 
Resources: \textbf{M.J.K.};
Funding acquisition: \textbf{M.J.K.}. 
All authors read and approved the final manuscript.


\begin{thebibliography}{10}

\bibitem{bentley1986experimental}
B.~J. Bentley and L.~G. Leal.
\newblock {An experimental investigation of drop deformation and breakup in
  steady, two-dimensional linear flows}.
\newblock {\em Journal of Fluid Mechanics}, 167:241--283, 1986.
\newblock \href {https://doi.org/10.1017/S0022112086002811}
  {\path{doi:10.1017/S0022112086002811}}.

\bibitem{bhatnagar1954model}
P.~L. Bhatnagar, E.~P. Gross, and M.~Krook.
\newblock {A Model for Collision Processes in Gases. I. Small Amplitude
  Processes in Charged and Neutral One-Component Systems}.
\newblock {\em Physical Review E}, 94:511--525, 1954.
\newblock \href {https://doi.org/10.1103/PhysRev.94.511}
  {\path{doi:10.1103/PhysRev.94.511}}.

\bibitem{bigio1998predicting}
D.~I. Bigio, C.~R. Marks, and R.~V. Calabrese.
\newblock {Predicting drop breakup in complex flows from model flow
  experiments}.
\newblock {\em International Polymer Processing}, 13(2):192--198, 1998.
\newblock \href {https://doi.org/10.3139/217.980192}
  {\path{doi:10.3139/217.980192}}.

\bibitem{bukreev2022consistent}
F.~Bukreev, S.~Simonis, A.~Kummerl{\"a}nder, J.~Je{\ss}berger, and M.~J.
  Krause.
\newblock {Consistent lattice Boltzmann methods for the volume averaged
  Navier\textendash Stokes equations}.
\newblock {\em arXiv}, preprint(arXiv:2208.09267), 2022.
\newblock \href {https://doi.org/10.48550/arXiv.2208.09267}
  {\path{doi:10.48550/arXiv.2208.09267}}.

\bibitem{calhoun2022systematic}
S.~G.~K. Calhoun, K.~K. Brower, V.~C. Suja, G.~Kim, N.~Wang, A.~L. McCully,
  H.~Kusumaatmaja, G.~G. Fuller, and P.~M. Fordyce.
\newblock {Systematic characterization of effect of flow rates and buffer
  compositions on double emulsion droplet volumes and stability}.
\newblock {\em Lab on a Chip}, 22:2315--2330, 2022.
\newblock \href {https://doi.org/10.1039/D2LC00229A}
  {\path{doi:10.1039/D2LC00229A}}.

\bibitem{dapelo2021lattice-boltzmann}
D.~Dapelo, S.~Simonis, M.~J. Krause, and J.~Bridgeman.
\newblock {Lattice-Boltzmann coupled models for advection–diffusion flow on a
  wide range of Péclet numbers}.
\newblock {\em Journal of Computational Science}, 51:101363, 2021.
\newblock \href {https://doi.org/10.1016/j.jocs.2021.101363}
  {\path{doi:10.1016/j.jocs.2021.101363}}.

\bibitem{guo2002forcing}
Z.~Guo, C.~Zheng, and B.~Shi.
\newblock {Discrete lattice effects on the forcing term in the lattice
  Boltzmann method}.
\newblock {\em Physical Review E}, 65:046308, 2002.
\newblock \href {https://doi.org/10.1103/PhysRevE.65.046308}
  {\path{doi:10.1103/PhysRevE.65.046308}}.

\bibitem{haussmann2021fluid-structure}
M.~Haussmann, P.~Reinshaus, S.~Simonis, H.~Nirschl, and M.~J. Krause.
\newblock {Fluid{\textendash}Structure Interaction Simulation of a Coriolis
  Mass Flowmeter Using a Lattice Boltzmann Method}.
\newblock {\em Fluids}, 6(4):167, 2021.
\newblock \href {https://doi.org/10.3390/fluids6040167}
  {\path{doi:10.3390/fluids6040167}}.

\bibitem{haussmann2019direct}
M.~Haussmann, S.~Simonis, H.~Nirschl, and M.~J. Krause.
\newblock Direct numerical simulation of decaying homogeneous isotropic
  turbulence -- numerical experiments on stability, consistency and accuracy of
  distinct lattice {Boltzmann methods}.
\newblock {\em International Journal of Modern Physics C}, 30(09):1--29, 2019.
\newblock \href {https://doi.org/10.1142/S0129183119500748}
  {\path{doi:10.1142/S0129183119500748}}.

\bibitem{higdon1993kinematics}
J.~J.~L. Higdon.
\newblock {The kinematics of the four-roll mill}.
\newblock {\em Physics of Fluids A}, 5(1):274--276, 1993.
\newblock \href {https://doi.org/10.1063/1.858782}
  {\path{doi:10.1063/1.858782}}.

\bibitem{park2019taylor}
V.~T. Hoang and J.~M. Park.
\newblock {A Taylor analogy model for droplet dynamics in planar extensional
  flow}.
\newblock {\em Chemical Engineering Science}, 204:27--34, 2019.
\newblock \href {https://doi.org/10.1016/j.ces.2019.04.015}
  {\path{doi:10.1016/j.ces.2019.04.015}}.

\bibitem{hsu2009deformation}
A.~S. Hsu and L.~G. Leal.
\newblock {Deformation of a viscoelastic drop in planar extensional flows of a
  Newtonian fluid}.
\newblock {\em Journal of Non-Newtonian Fluid Mechanics}, 160(2-3):176--180,
  2009.
\newblock \href {https://doi.org/10.1016/j.jnnfm.2009.03.004}
  {\path{doi:10.1016/j.jnnfm.2009.03.004}}.

\bibitem{huang2015multiphase}
H.~Huang, M.~Sukop, and X.~Lu.
\newblock {\em {Multiphase lattice Boltzmann methods: Theory and application}}.
\newblock John Wiley {\&} Sons, 2015.
\newblock \href {https://doi.org/10.1002/9781118971451}
  {\path{doi:10.1002/9781118971451}}.

\bibitem{kendon2001inertial}
V.~M. Kendon, M.~E. Cates, I.~Pagonabarraga, J.-C. Desplat, and P.~Bladon.
\newblock {Inertial effects in three-dimensional spinodal decomposition of a
  symmetric binary fluid mixture: a lattice Boltzmann study}.
\newblock {\em Journal of Fluid Mechanics}, 440:147--203, 2001.
\newblock \href {https://doi.org/10.1017/S0022112001004682}
  {\path{doi:10.1017/S0022112001004682}}.

\bibitem{kim2007phase}
J.~Kim.
\newblock {Phase field computations for ternary fluid flows}.
\newblock {\em Computer Methods in Applied Mechanics and Engineering},
  196(45-48):4779--4788, 2007.
\newblock \href {https://doi.org/10.1016/j.cma.2007.06.016}
  {\path{doi:10.1016/j.cma.2007.06.016}}.

\bibitem{komrakova2014lattice}
A.~E. Komrakova, O.~Shardt, D.~Eskin, and J.~J. Derksen.
\newblock {Lattice Boltzmann simulations of drop deformation and breakup in
  shear flow}.
\newblock {\em International Journal of Multiphase Flow}, 59:24--43, 2014.
\newblock \href {https://doi.org/10.1016/j.ijmultiphaseflow.2013.10.009}
  {\path{doi:10.1016/j.ijmultiphaseflow.2013.10.009}}.

\bibitem{krause2021openlb}
M.~J. Krause, A.~Kummerländer, S.~J. Avis, H.~Kusumaatmaja, D.~Dapelo,
  F.~Klemens, M.~Gaedtke, N.~Hafen, A.~Mink, R.~Trunk, J.~E. Marquardt, M.-L.
  Maier, M.~Haussmann, and S.~Simonis.
\newblock {OpenLB{\textemdash}Open source lattice Boltzmann code}.
\newblock {\em Computers {\&} Mathematics with Applications}, 81:258--288,
  2021.
\newblock \href {https://doi.org/10.1016/j.camwa.2020.04.033}
  {\path{doi:10.1016/j.camwa.2020.04.033}}.

\bibitem{krause2020openlb14}
M.J. Krause, S.~Avis, H.~Kusumaatmaja, D.~Dapelo, M.~Gaedtke, N.~Hafen,
  M.~Haußmann, Jonathan Jeppener-Haltenhoff, L.~Kronberg, A.~Kummerländer,
  J.E. Marquardt, T.~Pertzel, S.~Simonis, R.~Trunk, M.~Wu, and A.~Zarth.
\newblock {OpenLB Release 1.4: Open Source Lattice Boltzmann Code}, November
  2020.
\newblock \href {https://doi.org/10.5281/zenodo.4279263}
  {\path{doi:10.5281/zenodo.4279263}}.

\bibitem{kruger2017lattice}
T.~Kr{\"u}ger, H.~Kusumaatmaja, A.~Kuzmin, O.~Shardt, G.~Silva, and E.~M.
  Viggen.
\newblock {\em {The Lattice Boltzmann Method: Principles and Practice}}.
\newblock Springer International Publishing, 2017.
\newblock \href {https://doi.org/10.1007/978-3-319-44649-3}
  {\path{doi:10.1007/978-3-319-44649-3}}.

\bibitem{kummerlander2022implicit}
A.~Kummerl{\"a}nder, M.~Dorn, M.~Frank, and M.~J. Krause.
\newblock {Implicit Propagation of Directly Addressed Grids in Lattice
  Boltzmann Methods}.
\newblock {\em Concurrency and Computation}, accepted, 2021.
\newblock \href {https://doi.org/10.1002/cpe.7509}
  {\path{doi:10.1002/cpe.7509}}.

\bibitem{kusumaatmaja2010lattice}
H.~Kusumaatmaja and J.~M. Yeomans.
\newblock {Lattice Boltzmann simulations of wetting and drop dynamics}.
\newblock In {\em Simulating Complex Systems by Cellular Automata}, pages
  241--274. Springer, 2010.
\newblock \href {https://doi.org/10.1007/978-3-642-12203-3_11}
  {\path{doi:10.1007/978-3-642-12203-3_11}}.

\bibitem{lallemand2021lattice}
P.~Lallemand, L.-S. Luo, M.~Krafczyk, and W.-A. Yong.
\newblock {The lattice Boltzmann method for nearly incompressible flows}.
\newblock {\em Journal of Computational Physics}, 431:109713, 2021.
\newblock \href {https://doi.org/10.1016/j.jcp.2020.109713}
  {\path{doi:10.1016/j.jcp.2020.109713}}.

\bibitem{li2000numerical}
J.~Li, Y.~Y. Renardy, and M.~Renardy.
\newblock {Numerical simulation of breakup of a viscous drop in simple shear
  flow through a volume-of-fluid method}.
\newblock {\em Physics of Fluids}, 12(2):269--282, 2000.
\newblock \href {https://doi.org/10.1063/1.870305}
  {\path{doi:10.1063/1.870305}}.

\bibitem{mei2006consistent}
R.~Mei, L.-S. Luo, P.~Lallemand, and D.~d’Humi{\`e}res.
\newblock {Consistent initial conditions for lattice {Boltzmann} simulations}.
\newblock {\em Computers {\&} Fluids}, 35(8-9):855--862, 2006.
\newblock \href {https://doi.org/10.1016/j.compfluid.2005.08.008}
  {\path{doi:10.1016/j.compfluid.2005.08.008}}.

\bibitem{mink2021comprehensive}
A.~Mink, K.~Schediwy, C.~Posten, H.~Nirschl, S.~Simonis, and M.~J. Krause.
\newblock {Comprehensive Computational Model for Coupled Fluid Flow, Mass
  Transfer, and Light Supply in Tubular Photobioreactors Equipped with Glass
  Sponges}.
\newblock {\em Energies}, 15(20), 2022.
\newblock \href {https://doi.org/10.3390/en15207671}
  {\path{doi:10.3390/en15207671}}.

\bibitem{semprebon2016ternary}
C.~Semprebon, T.~Kr{\"u}ger, and H.~Kusumaatmaja.
\newblock {Ternary free-energy lattice Boltzmann model with tunable surface
  tensions and contact angles}.
\newblock {\em Physical Review E}, 93(3):033305, 2016.
\newblock \href {https://doi.org/10.1103/PhysRevE.93.033305}
  {\path{doi:10.1103/PhysRevE.93.033305}}.

\bibitem{shapira1990low}
M.~Shapira and S.~Haber.
\newblock {Low Reynolds number motion of a droplet in shear flow including wall
  effects}.
\newblock {\em International Journal of Multiphase Flow}, 16(2):305--321, 1990.
\newblock \href {https://doi.org/10.1016/0301-9322(90)90061-M}
  {\path{doi:10.1016/0301-9322(90)90061-M}}.

\bibitem{simonis2020relaxation}
S.~Simonis, M.~Frank, and M.~J. Krause.
\newblock {On relaxation systems and their relation to discrete velocity
  Boltzmann models for scalar advection--diffusion equations}.
\newblock {\em Philosophical Transactions of the Royal Society of London,
  Series A: Mathematical, Physical and Engineering Sciences},
  378(2175):20190400, 2020.
\newblock \href {https://doi.org/10.1098/rsta.2019.0400}
  {\path{doi:10.1098/rsta.2019.0400}}.

\bibitem{simonis2022constructing}
S.~Simonis, M.~Frank, and M.~J. Krause.
\newblock {Constructing relaxation systems for lattice Boltzmann methods}.
\newblock {\em {Applied Mathematics Letters}}, 137:108484, 2023.
\newblock \href {https://doi.org/10.1016/j.aml.2022.108484}
  {\path{doi:10.1016/j.aml.2022.108484}}.

\bibitem{simonis2021linear}
S.~Simonis, M.~Haussmann, L.~Kronberg, W.~D{\"o}rfler, and M.~J. Krause.
\newblock {Linear and brute force stability of orthogonal moment
  multiple-relaxation-time lattice Boltzmann methods applied to homogeneous
  isotropic turbulence}.
\newblock {\em Philosophical Transactions of the Royal Society of London,
  Series A: Mathematical, Physical and Engineering Sciences},
  379(2208):20200405, 2021.
\newblock \href {https://doi.org/10.1098/rsta.2020.0405}
  {\path{doi:10.1098/rsta.2020.0405}}.

\bibitem{simonis2022forschungsnahe}
S.~Simonis and M.~J. Krause.
\newblock {Forschungsnahe Lehre unter Pandemiebedingungen}.
\newblock {\em Mitteilungen der Deutschen Mathematiker-Vereinigung},
  30(1):43--45, 2022.
\newblock \href {https://doi.org/10.1515/dmvm-2022-0015}
  {\path{doi:10.1515/dmvm-2022-0015}}.

\bibitem{simonis2022limit}
S.~Simonis and M.~J. Krause.
\newblock {Limit Consistency for Lattice Boltzmann Equations}.
\newblock {\em arXiv}, preprint(arXiv:2208.06867), 2022.
\newblock \href {https://doi.org/10.48550/arXiv.2208.06867}
  {\path{doi:10.48550/arXiv.2208.06867}}.

\bibitem{simonis2022temporal}
S.~Simonis, D.~Oberle, M.~Gaedtke, P.~Jenny, and M.~J. Krause.
\newblock {Temporal large eddy simulation with lattice Boltzmann methods}.
\newblock {\em Journal of Computational Physics}, 454:110991, 2022.
\newblock \href {https://doi.org/10.1016/j.jcp.2022.110991}
  {\path{doi:10.1016/j.jcp.2022.110991}}.

\bibitem{siodlaczek2021numerical}
M.~Siodlaczek, M.~Gaedtke, S.~Simonis, M.~Schweiker, N.~Homma, and M.~J.
  Krause.
\newblock {Numerical evaluation of thermal comfort using a large eddy lattice
  Boltzmann method}.
\newblock {\em Building and Environment}, 192:107618, 2021.
\newblock \href {https://doi.org/10.1016/j.buildenv.2021.107618}
  {\path{doi:10.1016/j.buildenv.2021.107618}}.

\bibitem{soligo2020deformation}
G.~Soligo, A.~Roccon, and A.~Soldati.
\newblock {Deformation of clean and surfactant-laden droplets in shear flow}.
\newblock {\em Meccanica}, 55(2):371--386, 2020.
\newblock \href {https://doi.org/10.1007/s11012-019-00990-9}
  {\path{doi:10.1007/s11012-019-00990-9}}.

\bibitem{swift1996lattice}
M.~R. Swift, E.~Orlandini, W.~R. Osborn, and J.~M. Yeomans.
\newblock {Lattice Boltzmann simulations of liquid-gas and binary fluid
  systems}.
\newblock {\em Physical Review E}, 54(5):5041, 1996.
\newblock \href {https://doi.org/10.1103/PhysRevE.54.5041}
  {\path{doi:10.1103/PhysRevE.54.5041}}.

\bibitem{taylor1934formation}
G.~I. Taylor.
\newblock {The formation of emulsions in definable fields of flow}.
\newblock {\em Proceedings of the Royal Society of London, Series A:
  Mathematical and Physical Sciences}, 146(858):501--523, 1934.
\newblock \href {https://doi.org/10.1098/rspa.1934.0169}
  {\path{doi:10.1098/rspa.1934.0169}}.

\bibitem{tretheway2001deformation}
D.~C. Tretheway and L.~G. Leal.
\newblock {Deformation and relaxation of Newtonian drops in planar extensional
  flows of a Boger fluid}.
\newblock {\em Journal of Non-Newtonian Fluid Mechanics}, 99(2-3):81--108,
  2001.
\newblock \href {https://doi.org/10.1016/S0377-0257(01)00123-9}
  {\path{doi:10.1016/S0377-0257(01)00123-9}}.

\bibitem{wang2020modelling}
N.~Wang, Ciro Semprebon, Haihu Liu, Chuhua Zhang, and Halim Kusumaatmaja.
\newblock {Modelling double emulsion formation in planar flow-focusing
  microchannels}.
\newblock {\em Journal of Fluid Mechanics}, 895:A22, 2020.
\newblock \href {https://doi.org/10.1017/jfm.2020.299}
  {\path{doi:10.1017/jfm.2020.299}}.

\bibitem{zhao2007drop}
X.~Zhao.
\newblock {Drop breakup in dilute Newtonian emulsions in simple shear flow: New
  drop breakup mechanisms}.
\newblock {\em Journal of Rheology}, 51(3):367--392, 2007.
\newblock \href {https://doi.org/10.1122/1.2714641}
  {\path{doi:10.1122/1.2714641}}.

\bibitem{zheng2014shrinkage}
L.~Zheng, T.~Lee, Z.~Guo, and D.~Rumschitzki.
\newblock {Shrinkage of bubbles and drops in the lattice Boltzmann equation
  method for nonideal gases}.
\newblock {\em Physical Review E}, 89:033302, 2014.
\newblock \href {https://doi.org/10.1103/PhysRevE.89.033302}
  {\path{doi:10.1103/PhysRevE.89.033302}}.

\end{thebibliography}
\end{document}